\documentclass[a4paper,11pt]{article}
\usepackage{graphicx}
\usepackage{amssymb}
\usepackage{latexsym, bm}
\usepackage{multicol}
\usepackage{indentfirst}
\usepackage{amssymb,amsfonts}
\usepackage{amsmath}
\usepackage{setspace}
\usepackage{enumerate}
\usepackage{CJK}
\newtheorem{theorem}{Theorem}[section]
\newtheorem{lemma}{Lemma}[section]
\newtheorem{claim}{Claim}[section]
\newtheorem{fact}{Fact}[section]
\newtheorem{corollary}{Corollary}[section]
\newtheorem{prop}{Proposition}[section]
\newtheorem{definition}{Definition}[section]

\usepackage{amssymb}
\begin{spacing}{1.0}
\begin{document}
\title{Three-color Ramsey number of an odd cycle versus bipartite graphs
with small bandwidth \footnote{Center for Discrete Mathematics, Fuzhou University, Fuzhou, 350108, P.~R.~China.  Email: {\em chunlin\_you@163.com, linqizhong@fzu.edu.cn. }  }}
\date{}
\author{Chunlin You \;\; and \;\; Qizhong Lin\footnote{Corresponding author. Supported in part by NSFC (No.\ 12171088).}
}

\maketitle
\begin{abstract}

A graph $\mathcal{H}=(W,E_\mathcal{H})$ is said to have {\em bandwidth} at most $b$ if there exists a labeling of $W$ as $w_1,w_2,\dots,w_n$ such that $|i-j|\leq b$ for every edge $w_iw_j\in E_\mathcal{H}$.
We say that $\mathcal{H}$ is a {\em balanced  $(\beta,\Delta)$-graph} if it is a bipartite graph with bandwidth at most $\beta |W|$ and maximum degree at most $\Delta$, and it also has a proper 2-coloring $\chi :W\rightarrow[2]$ such that  $||\chi^{-1}(1)|-|\chi^{-1}(2)||\leq\beta|\chi^{-1}(2)|$.

In this paper, we prove that for every $\gamma>0$ and every natural number $\Delta$,
there exists a constant $\beta>0$ such that for every
balanced   $(\beta,\Delta)$-graph $\mathcal{H}$ on $n$ vertices we have  $$R(\mathcal{H},  \mathcal{H}, C_n) \leq (3+\gamma)n$$ for all sufficiently large odd $n$. The upper bound is sharp for several classes of graphs.
Let $\theta_{n,t}$ be the graph consisting of $t$ internally disjoint paths of length $n$ all sharing the same endpoints. As a corollary, for each fixed $t\geq 1$, $R(\theta_{n, t},\theta_{n, t}, C_{nt+\lambda})=(3t+o(1))n,$ where $\lambda=0$ if $nt$ is odd and $\lambda=1$ if $nt$ is even. In particular, we have $R(C_{2n},C_{2n}, C_{2n+1})=(6+o(1))n$, which is a special case of a result of Figaj and {\L}uczak (2018).

\medskip

{\bf Keywords:} \  Ramsey number;   Small bandwidth; Cycle; Regularity Lemma

\medskip

{\bf Mathematics Subject Classification:} \  05C55;\;\;05D10

\end{abstract}

\section{Introduction}
%Throughout this paper, we only concern with undirected simple finite graphs.
For graphs $G$ and $G_1,G_2,\dots,G_k$ we write $G\rightarrow(G_1,G_2,\dots,G_k)$ if
for every $k$-edge-coloring of $G$, there is a monochromatic copy of $G_i$ for some $i\in[k]$, where $[k]=\{1,2,\dots,k\}$.
The multicolor Ramsey number $R(G_1,G_2,\dots,G_k)$ is defined to be the smallest integer $N$ such that $K_N\rightarrow(G_1,G_2,\dots,G_k)$.
We denote $R_k(G)$ instead of $R(G_1,G_2,\dots,G_k)$ when $G_i=G$ for $1\le i\le k$.

Let $P_n$ (resp. $C_n$) be the path (resp. cycle) on $n$ vertices.
For the Ramsey numbers of $R(C_m,C_n)$, it has been studied and completely determined in Bondy and Erd\H{o}s \cite{b-e}, Faudree and Schelp \cite{f-s}, and Rosta \cite{ros}. Later, Erd\H{o}s, Faudree, Rousseau and  Schelp \cite{efrs} obtained the exact value of $R(C_\ell,C_m,C_n)$ for sufficiently large $n$. In 1999, by using the regularity lemma, {\L}uczak \cite{lucazk-1999} proved that $R(C_n,C_n,C_n)=(4+o(1))n$ for all large odd integer $n$, and some years later Kohayakawa, Simonovits and Skokan \cite{kss-2005} determined the exact value of it for all large odd $n$.
In \cite{lucazk-1999}, {\L}uczak introduced a technique that uses the  regularity lemma to reduce
problems about paths and cycles to problems about
\emph{connected matchings}, which are
matchings that are contained in a connected component.
This technique has become fairly standard.
The following papers were based on this technique.
For large even cycle, Figaj and {\L}uczak \cite{f-l-2007} obtained that $R(C_n,C_n,C_n)=(2+o(1))n$, and  Benevides and Skokan \cite{Benevides-Skokan-2009} obtained the exact value. Furthermore, Figaj and  {\L}uczak \cite{f-l-2007,Figaj-luczak-2018} determined the asymptotic values of $R(C_\ell,C_m,C_n)$ when $\ell,m$, and $n$ are large integers with same order. %Recently, Jenssen and Skokan \cite{j-s} established that $R_k(C_n)=2^{k-1}(n-1)+1$ for each fixed $k\ge2$ and large odd $n$, which confirms a conjecture of \cite{b-e}.
%In fact, Ferguson \cite{Ferguson} has proved a sharpening of the Figaj and {\L}uczak result, obtaining exact Ramsey numbers for the case of mixed parity cycles.
In particular, Figaj and {\L}uczak \cite{Figaj-luczak-2018}
determined \[R(C_{2n}, C_{2n}, C_{2n+1})=(6+o(1))n.\]
Furthermore, Ferguson \cite{Ferguson-1, Ferguson-2, Ferguson-3} obtained exact Ramsey numbers for the  three-color cases of mixed parity cycles.
Recently, Jenssen and Skokan \cite{j-s} established that
$$R_k(C_n)=2^{k-1}(n-1)+1$$ for each fixed $k\ge2$ and large odd $n$, which confirms a conjecture of Bondy and Erd\H{o}s \cite{b-e}.
For more multicolor Ramsey numbers involving large cycles, we refer the reader to \cite{ a-b-s-2013,efrs, lucss, sark-2016}, etc.

The Ramsey numbers of paths are also fruitful.
A well-known result of Gerencs\'{e}r and Gy\'{a}rf\'{a}s \cite{gerence-gyarfas-1967} states that
$R({P_m},{P_n}) = n+\lfloor {\frac{m}{2}}\rfloor$ for all $n\ge m\ge1$.
Faudree and Schelp \cite{Faudree-Schelp-1975} determined $R({P_{{n_0}}},{P_{n_1}},P_{n_2})$
for  $n_0\ge6(n_1+n_2)^2$, and they also conjectured that
 \begin{align}\label{path-path-p}
R({P_n},{P_n},{P_n}) =\left\{ {\begin{array}{*{20}{c}}
{2n- 1,}\\{2n- 2,}
\end{array}} \right.\begin{array}{*{20}{c}}
\text{$n$ is odd},\\
\text{$n$ is even}.
\end{array}
\end{align}
Gy\'{a}rf\'{a}s, Ruszink\'{o}, S\'{a}rk\"{o}zy and Szemer\'{e}di \cite{gyarfas-szemerdi-2007} confirmed this conjecture for large $n$.
For Ramsey numbers of cycles versus paths, Faudree, Lawrence, Parsons and Schelp \cite{FLPS}  obtained the Ramsey numbers for all path-cycle pairs.
%For multicolor Ramsey numbers involving large cycles, see also \cite{efrs}.
In  2009,  Dzido and Fidytek \cite{DzFi2} (independently by Bielak \cite{Bielak}) obtained that if $m\ge3$ is odd, $n\geq m$ and $n > \frac{{3{t^2} - 14t + 25}}{4}$
when $t$ is odd, and $n > \frac{{3{t^2} - 10t + 20}}{8}$ when $t$ is even, then
$R({P_t},{P_n},{C_m})= 2n + 2\left\lfloor {\frac{t}{2}} \right\rfloor  - 3.$
In \cite{shao-xu}, Shao, Xu, Shi and Pan proved that $R(P_4, P_5, C_n)=n+2$ for $n\ge23$ and $R(P_4, P_6, C_n)=n+3$  for $n\ge18$. Omidi and Raeisi \cite{Omidi-Raeisi-2011} determined $R(P_4,P_n,C_4)=R(P_5,P_n,C_4)=n+2$ for $n\ge5$.

%For more cycle-path Ramsey numbers, see, e.g.,  \cite{DzNS, lucss,Omidi-Raeisi-2011, sark-2016, Davies-2017}.

We say that  a graph $\mathcal{H}=(W,E_\mathcal{H})$ has bandwidth at most $b$ if there is a labeling of $W$ as $w_1,w_2,\dots,w_n$ such that $|i-j|\leq b$ for every edge $w_iw_j\in E_\mathcal{H}$ (see, e.g., \cite{julia-klaas}).
We say $\mathcal{H}=(W,E_\mathcal{H})$ is a {\em balanced   $(\beta,\Delta)$-graph} if it is a  bipartite graph with bandwidth at most $\beta |W|$ and maximum degree at most $\Delta$, and furthermore it has a proper 2-coloring
$\chi :W\rightarrow[2]$ such that  $||\chi^{-1}(1)|-|\chi^{-1}(2)||\leq\beta|\chi^{-1}(2)|$.
%For example, it was shown in \cite{julia-klaas}  that sufficiently large planar graphs with maximum degree at most $\Delta$ are $(\beta,\Delta)$-graphs  for any fixed $\beta>0$.
% Mota et al. \cite{G-G-M} determined asymptotically the two color Ramsey number for bipartite graphs with small bandwidth and bounded maximum degree and the three color Ramsey number for balanced $(\beta,\Delta)$-graphs   $\mathcal{H}$ was estimated in \cite{mota,G-G-M}.
 In \cite{G-G-M}, Mota, S\'{a}rk\"{o}zy, Schacht and  Taraz determined asymptotically the three color  Ramsey number $R(\mathcal{H},\mathcal{H},\mathcal{H})\leq (2+o(1))n$, which implies that $$R(C_{n},C_n,C_n)=(2+o(1))n$$ for all large even $n$.
%for every  balanced   $(\beta,\Delta)$-graph $\mathcal{H}$ on $n$ vertices.
%Later, Shen, Lin and Liu \cite{shen-lin-liu} obtained the bipartite Ramsey numbers for bipartite graphs of small bandwidth.
%Denote $H_{n-1}$ be the graph $H$ on $n-1$ vertices.

In this paper, we are concerned with the asymptotic behavior of the Ramsey number $R(\mathcal{H},\mathcal{H},C_n)$ when  $n$  is odd and every  balanced $(\beta,\Delta)$-graph $\mathcal{H}$ on  $n$ vertices.

\begin{theorem}\label{main theorem}
For any sufficiently small $\gamma$ and every natural number $\Delta$,
there exists a constant $\beta>0$ such that for every balanced   $(\beta,\Delta)$-graph $\mathcal{H}$ on $n$ vertices we have  $$R(\mathcal{H},  \mathcal{H}, C_n) \leq (3+\gamma)n$$ for all sufficiently large odd $n$.
%The inequality is asymptotic tight if $\mathcal{H}$ is a connected balanced $(\beta,\Delta)$-graph on $n$ vertices.
\end{theorem}

%{\em Remark.} If $\mathcal{H}$ is a connected balanced $(\beta,\Delta)$-graph on $n$ vertices, then $\mathcal{H}$ contains a spanning tree as a subgraph. Thus we have $R(T_n,T_n, C_n)>3n-5$ for all odd $n$. This implies that the upper bound in Theorem \ref{main theorem} is asymptotic tight if $\mathcal{H}$ is a connected balanced $(\beta,\Delta)$-graph on $n$ vertices.

\medskip
 Let $\theta_{n,t}$ be the graph consisting of $t$ internally disjoint paths of length $n$ all sharing the same endpoints. In particular, $\theta_{n,2}=C_{2n}$, i.e. the even cycle on $2n$ vertices.
Observe that $\theta_{n,t}$ is a bipartite graph with $(n-1)t +2$ vertices,
and $\theta_{n,t}$ is a  balanced $(\beta,t)$-graph with $\beta \leq\frac{2}{n}$.

%Denote $m\sim n$ if $\lim\limits_{n\to \infty}\frac{m}{n}=1$.
The following result is clear from Theorem \ref{main theorem}.
\begin{corollary}\label{corollary1}
For each fixed $t\geq 1$,
$$ R(\theta_{n, t},\theta_{n, t}, C_{nt+\lambda})=(3t+o(1))n,$$
where $\lambda=0$ if $nt$ is odd and $\lambda=1$ if $nt$ is even.
In particular, we have that $R(P_{2n},P_{2n}, C_{2n+1})\sim R(C_{2n},C_{2n}, C_{2n+1})=(6+o(1))n$.
\end{corollary}

Let us point out that the result $R(C_{2n},C_{2n}, C_{2n+1})=(6+o(1))n$
and $R(P_{2n},P_{2n}, C_{2n+1})=(6+o(1))n$
follow from the work of Figaj and {\L}uczak \cite{Figaj-luczak-2018}
and have been strengthened by Ferguson \cite{Ferguson-1}.

\medskip

\section{Preliminaries}\label{chap2}
In this paper, we shall omit all floors and ceilings since it will not affect our argument.
For a graph $G$, let $V(G)$ and $E(G)$ denote its vertex set and edge set, respectively.
Let $|E(G)|=e(G)$ and $|V(G)|=v(G)$.
For $v\in V$, we let $N_G(v)$ denote the neighborhood of $v$ in $G$,
and $\deg_G(v)=|N_G(v)|$ which is the degree of a vertex $v\in V(G)$.
We denote by $\delta(G)$  and $\triangle(G)$ the minimum and maximum
degrees of the vertices of $G$.
For any subset $U \subseteq V$, we use $G[U]$ to denote the subgraph induced by the vertex set $U$ in $G$.
For a vertex $v\in V$ and $U\subset V$, we write $N_G(v,U)$ for the neighbors of $v$ in $U$ in graph $G$ and $\deg(v,U)=|N_G(v,U)|$.
%For disjoint vertex sets $A, B\subseteq V$, let $e_G(A,B)$ denote the number of edges of $G$ with one endpoint in $A$ and the other in $B$, and the density between $A$ and $B$ is
%\[d_G(A,B)=\dfrac{e_G(A,B)}{|A| |B|}.\]
%We always omit the subscript when there is no confusion.

\subsection{The regularity method}
For disjoint vertex sets $A, B\subseteq V$, let $e_G(A,B)$ denote the number of edges of $G$ with one endpoint in $A$ and the other in $B$, and the density between $A$ and $B$ is
\[d_G(A,B)=\dfrac{e_G(A,B)}{|A| |B|}.\]
We always omit the subscript when there is no confusion.

\begin{definition}[$\epsilon$-regular]\label{regular}
A pair $(A,B)$ is $\epsilon$-regular if for all $X\subseteq A$ and $Y\subseteq B$ with $|X|>\epsilon|A|$ and $|Y|>\epsilon|B|$ we have $|d(X,Y)-d(A,B)|<\epsilon$.
\end{definition}
\begin{definition}[($\epsilon,d$)-regular]
For $\epsilon>0, d\leq 1$,  a pair $(A,B)$ is said to be $(\epsilon,d)$-regular if it is $\epsilon$-regular and $d(A,B)\geq d$.
\end{definition}
\begin{definition}[$(\epsilon,d)$-super-regular]
A pair $(A,B)$ is said to be $(\epsilon,d)$-super-regular if it is $\epsilon$-regular and $\deg(u,B)>d|B|$ for all $u\in A$ and $\deg(v,A)>d|A|$ for all $v\in B$.
\end{definition}

We have a property that any regular pair has a large subgraph which is super-regular,
and we include a proof for completeness.
\begin{fact}\label{localregular}
For $0<\epsilon<1/2$ and $d\leq1$, if $(A,B)$ is $(\epsilon,d)$-regular with $|A|=|B|=m$, then there exist $A_1\subseteq A$ and $B_1\subseteq B$ with $|A_1|=|B_1|\ge(1-\epsilon)m$ such that $(A_1,B_1)$ is $(2\epsilon,d-2\epsilon)$-super-regular.
\end{fact}

\noindent{\bf Proof.}
Let $X\subseteq A$ consist of vertices with at most $(d-\epsilon)|B|$ neighbors in $B$. Accounting for $e(X,B)\leq |X|\cdot(d-\epsilon)|B|$, we have  $|d(X,B)-d|\geq \epsilon$.
Due to the definition of $\epsilon$-regular, it follows that $|X|\le\epsilon m$.
Similarly, assume that $Y\subseteq B$ consists of vertices with at most $(d-\epsilon)|A|$ neighbors in $A$, then we obtain $|Y|\le\epsilon m$.
Take $A_1\subseteq A\setminus X$ and $B_1\subseteq B\setminus Y$ with $|A_1|=|B_1|=(1-\epsilon)m$.
Thus for any vertex $u\in A_1$, $\deg(u,B_1)\geq (d-2\epsilon)m$.
Also, $\deg(w,A_1)\geq (d-2\epsilon)m$ for any $w\in B_1$.

Now, for any subsets $S\subseteq A_1$ and $T\subseteq B_1$,
if $|S|> 2\epsilon|A_1|$ and $|T|> 2\epsilon|B_1|$,
then clearly $|S|> \epsilon m$ and $|T|> \epsilon m$.
Since $(A,B)$ is $(\epsilon,d)$-regular, it follows that
\[
|d(S,T)-d(A_1,B_1)|\le|d(S,T)-d(A,B)|+|d(A_1,B_1)-d(A,B)|<2\epsilon.
\]
We have thus proved the fact.
\hfill$\Box$

\medskip

In this paper, we will use the following  three-color version of the regularity lemma.  For many applications, we refer the reader to surveys \cite{ks,rs} and many recent references \cite{afz,cly,conlon,cf,cfw,flz, lp,nr09,sz15}, etc.

\begin{lemma}[Szemer\'{e}di \cite{regular-lemma}]\label{regular lemma}
 For every $\epsilon> 0$ and integer $t_0\geq 1$, there exists $T_0=T_0(\epsilon,t_0)\geq t_0$
such that the following holds. For all graphs $G_1, G_2$ and $G_3$ with the same vertex set $V$ and $|V|\geq t_0$, there exists a partition $V = \cup_{i=0}^t{V_i}$ satisfying $t_0\le t\le T_0$ and
\smallskip

$(1)$ $\left| {{V_0}} \right| < \epsilon n$, $\left| {{V_1}} \right| = \left| {{V_2}} \right| =  \ldots  = \left| {{V_t}} \right|$;

\smallskip

$(2)$  all but at most $\epsilon {t\choose 2}$ pairs $(V_i,V_j)$, $1 \le i \neq j \le t$, are $\epsilon$-regular for $G_1, G_2$ and $G_3$.
\end{lemma}

The next lemma by Benevides and Skokan \cite{Benevides-Skokan-2009} is a slightly stronger version compared to the original one established by {\L}uczak \cite[Claim 3]{lucazk-1999}.

\begin{lemma}[Benevides and Skokan \cite{Benevides-Skokan-2009}]\label{long-path-lemma}
For every $0<\beta_0<1$, there exists an $n_0$ such that for every $n>n_0$ the following holds: If $(V_1,V_2)$ is $\epsilon$-regular with $|V_1|=|V_2|=n$ and density at least $\beta_0/4$  for some $\epsilon$ satisfying $0<\epsilon<\beta_0/100$, then for every $\ell, 1\leq \ell \leq n-5\epsilon n/\beta_0$, and for every pair of vertices $v'\in V_1$, $v''\in V_2$ satisfying $\deg(v',V_2)$, $\deg(v'', V_1)\geq \beta_0 n/5$, $G$ contains a path of length $2\ell+1$ connecting $v'$ and $v''$.
\end{lemma}

Let $\mathcal{H}=(W,E_\mathcal{H})$ be a graph, for $S\subseteq W$, denote $N_\mathcal{H}(S)=[\cup_{v\in S} N_\mathcal{H}(v)]\setminus S$.
 For a graph $G = (V, E)$, a partition $\cup_{i = 1}^k {{V_i}}$
of $V$ is said to be $(\epsilon,d)$-regular on a reduced graph $H$ with vertex set contained in $[k]$ if the pair $(V_i,V_j)$ is $(\epsilon,d)$-regular whenever $ij\in E(H)$.

When applying the regularity lemma, we will indeed find
a partition of a monochromatic subgraph $G$ of $K_N$ with corresponding reduced graph containing a tree $T$ that contains a ``large'' matching $M$,
where the bipartite subgraphs of $G$ corresponding to the matching are super-regular pairs.
The following  definition of \emph{$\epsilon$-compatible} is due to Mota, S\'{a}rk\"{o}zy, Schacht and  Taraz \cite{G-G-M}.

\begin{definition}[$\epsilon$-compatible]\label{def}
Let $\mathcal{H}=(W, E_\mathcal{H})$ and $T=([t],E_{T})$ be graphs.
Let $M=([t],E_M)$ be a subgraph of $T$ where $E_M$ is a matching.
Given a partition  $W=\cup_{i = 1}^t {{W_i}}$, let $U_i$, for $i\in [t]$,
be the set of vertices in $W_i$, with neighbors in some $W_j$ with
$ij\in E_{T}\setminus E_{M}$.
 Set $U=\cup U_i$ and $U'_i=N_{\mathcal{H}}(U) \cap(W_{i}\setminus U)$.

We say that $W=\cup_{i = 1}^t {{W_i}}$ is $(\epsilon,T, M)$-compatible with a vertex partition $\cup_{i = 1}^t {{V_i}}$ of a graph $G=(V, E)$ if the following holds.

(1) $|W_i|\leq|V_i|$ for $i\in[t]$.

(2) $xy\in E_\mathcal{H}$ for $x\in W_i, y\in W_j$ implies $ij\in E_T$ for all distinct $i, j\in [t]$.

(3) $|U_i|\leq \epsilon |V_i|$ for $i\in[t]$.

(4) $|U'_i|, |U'_j|\leq \epsilon \min\{|V_i|, |V_j|: ij\in E_{M}\}$.
\end{definition}

The following corollary of the Blow-up Lemma (see B\"{o}ttcher, Heinig and Taraz \cite{J-B,J-P}) asserts that in the setup of Definition \ref{def} graphs $\mathcal{H}$ of bounded degree can be embedded into
$G$, if $G$ admits a partition being sufficiently regular on $T$ and super-regular on $M$.

\begin{lemma}[Embedding Lemma \cite{J-B,J-P}]\label{Embedding Lemma}
For all $d, \Delta>0$  there is a constant $\epsilon=\epsilon(d,\Delta)>0$ such that the following holds.
Let $G=(V,E)$ be an $N$-vertex graph that has a partition $\cup_{i = 1}^t{{V_i}}$ with $(\epsilon,d)$-reduced graph $T$ on $[t]$ which is $(\epsilon,d)$-super-regular on a graph $M\subset T$.
 Further, let  $\mathcal{H}=(W,E_\mathcal{H})$ be an
$n$-vertex graph with maximum degree $\Delta(\mathcal{H})\leq\Delta$ and
$n\leq N$ that has a vertex partition $\cup_{i = 1}^t{{W_i}}$ of $W$
 which is $(\epsilon, T, M)$-compatible with $\cup_{i = 1}^t{{V_i}}$.
Then $\mathcal{H}\subseteq G$.
\end{lemma}

For a graph $\mathcal{H}=(W,E_\mathcal{H})$ with $W=\{w_1,w_2,\dots,w_n\}$, where $w_i$ is a labeling of the vertices, let $\chi:W\rightarrow[2]$ be a 2-coloring.
For $W'\subseteq W$, denote $C_i(W')=|\chi^{-1}(i)\cap W'|$  for $i=1,2$. We know that $\chi$ is a $\beta$-balanced coloring of $W$ if $1-\beta\leq\frac{C_1(W)}{C_2(W)}\leq 1+\beta$.
%In this case, we say that $\mathcal{H}$  is globally balanced.
 A set $I\subseteq W$ is called an interval if there exists $p<q$ such that $I=\{w_p,w_{p+1},\dots,w_q\}$.
 Finally, let $\sigma:[\ell]\rightarrow[\ell]$ be a  permutation, and for a partition $\{I_1,I_2,\dots,I_\ell\}$ of $W$, where  $I_1,\dots,I_\ell$ are intervals,
let $C_i(\sigma,a,b)=\sum_{j=a}^bC_i(I_{\sigma(j)})$ for $i=1,2$.
%The following result  guarantees the local balancedness that we need.

\begin{lemma}[Mota et al. \cite{G-G-M}]\label{balanced}
For every $\xi>0$ and every integer $\ell\geq 1$ there exists $n_0$ such that if $\mathcal{H}=(W,E_\mathcal{H})$ is a graph on $W=\{w_1,w_2,\dots,w_n\}$ with $n\geq n_0$, then for every $\beta$-balanced 2-coloring $\chi$ of $W$ with $\beta\leq 2/\ell$, and every partition of $W$ into intervals $I_1,I_2,\dots,I_{\ell}$ with $|I_1|\leq|I_2|\leq\dots\leq|I_{\ell}|\leq|I_1|+1$ there exists a permutation $\sigma:[\ell]\rightarrow[\ell]$ such that for every pair of integers $1\leq a<b\leq \ell$ with $b-a\geq7/\xi$, we have
$|C_1(\sigma,a,b)-C_2(\sigma,a,b)|\leq\xi C_2(\sigma,a,b).$
\end{lemma}

\subsection{ Structure}
Recall that a set $M$ of independent edges in a graph $G=(V,E)$ is called a matching.
The size of a connected matching corresponds to the
number of edges in the matching.
The structure of graphs without large connected matchings play an important role in this paper.
As a path on $n$ vertices contains a connected matching on
$\lfloor n/2\rfloor$ edges, extremal
results for paths directly give an upper bound for connected matchings.

The following is the well-known extremal result for paths.
\begin{lemma}[Erd\H{o}s and Gallai \cite{Erd-Gallai-1959}]\label{gallai}
If $G$ is a graph of order $N$ which contains no $P_n$, then
\[e(G)\le \frac{{n - 2}}{2}N.\]
\end{lemma}

%\begin{corollary}\label{gallai}
%Let $G$ be a graph which does not contain a connected matching of size $n/2$ for even $n$. Then \[e(G)\le \frac{{n - 2}}{2}v(G).\]
%\end{corollary}

%\medskip

The following result describes the structure of a graph without a large matching.
\begin{lemma}[Knierim and Su \cite{kn-su}]\label{su-lemma}
For every connected graph $G = (V, E)$ which contains no matching of size $n/2$, there is a partition $S_G\cup Q_G \cup I_G$ of the vertex set $V$ such that
\smallskip

$(i)$ $\left| {{Q_G}} \right| + 2\left| {{S_G}} \right| = \min \left\{ {v(G),n - 1} \right\}$,

\smallskip

$(ii)$ $I_G$ is an independent set; additionally, if $v(G)\leq n-1$, then $I_G=\emptyset$,
	
\smallskip	

 $(iii)$ every vertex in $Q_G$ has at most one neighbor in $I_G$,

\smallskip

$(iv)$  every vertex in $I_G$ has degree less than $n/2$.
\end{lemma}

We will also use the following lemmas.
Let us begin with a result due to {\L}uczak \cite{lucazk-1999} which gives a description of the structure of a graph that contains no large odd cycle as a subgraph.
\begin{lemma}[{\L}uczak \cite{lucazk-1999}]\label{luc-3}
For every  $0 < \delta  < {10^{ - 15}}$, $\alpha\geq 2\delta$ and
$t\geq exp(\delta^{-16}/\alpha)$ the following holds. Each graph $H$ on $t$ vertices  which contains no odd cycles longer than $\alpha t$ contains subgraphs $H'$ and $H''$ such that:

\smallskip
$(i)$ $V(H') \cup V(H'')=V(H)$, $V(H')\cap V(H'')=\emptyset$ and each of the sets
$V(H')$ and $V(H'')$  is either empty or contains at least $\alpha \delta t/2$ vertices;

\smallskip
$(ii)$ $H'$ is bipartite;

\smallskip
$(iii)$ $H''$ contains no more than $\alpha t\left| {V(H'')} \right|/2$ edges;

\smallskip
$(iv)$  all except no more than $\delta {t^2}$ edges of $H$ belong to either $H'$ or $H''$.
\end{lemma}

The following result due to Gy\'{a}rf\'{a}s, Ruszink\'{o}, S\'{a}rk\"{o}zy and Szemer\'{e}di \cite{gyarfas-szemerdi-2007} states that a graph with high density always contains a dense subgraph.

\begin{fact}[Gy\'{a}rf\'{a}s et al. \cite{gyarfas-szemerdi-2007}]\label{density-degree}
Let $\epsilon>0$ be  sufficiently small
and let $H$ be a graph with $v(H)$ vertices.
If $e(H)\geq {v(H)\choose{2}}-\epsilon{t\choose{2}}$,  then $H$ has a subgraph $H'$ with at least $v(H) - \sqrt \epsilon  t$ vertices and $\delta(H')\geq v(H) - 2\sqrt \epsilon  t$.
\end{fact}

A spanning subgraph of a graph $G$ is a subgraph obtained by edge deletions only,
in other words, a subgraph whose vertex set is the entire vertex set of $G$.
We also need the following simple result.
%A similar argument as Fact \ref{density-degree} yields the following Lemma.
\begin{lemma}\label{density-degree-upper}
For any $\epsilon>0$, if $H'$ is a spanning subgraph of $ H$ with $e(H')\geq e(H)-\epsilon{t\choose{2}}$, then there exists a induced subgraph $H'' \subset H$ with at least
$v(H)- \sqrt \epsilon  t$ vertices and $\deg_{H'' }(u) <\deg_{H'}(u)+\sqrt \epsilon t$
for any vertex $u\in V(H'')$.
\end{lemma}
{\bf Proof.} Let $X=\{ u\in V( H)|  \deg_H(u)-\deg_{H'}(u)\geq \sqrt{\epsilon}t\}$.
Clearly, $e(H)-e(H')\geq\frac{\sqrt{\epsilon}t|X|}{2}$.
Thus we have \[\frac{\sqrt{\epsilon}t|X|}{2}\leq \epsilon {t\choose{2}},\]
implying that $|X|<\sqrt{\epsilon}t$. Denote $H''=H- X$, i.e. the subgraph obtained from $H$ by deleting all vertices of $X$ and all edges incident to some vertices of $X$.
Thus for any $u\in V(H'')$,
\[\deg_{H''}(u)=\deg_H(u)-\deg_H(u,X)<\left(\deg_{H'}(u)+\sqrt\epsilon  t \right)-\deg_H(u, X)\leq \deg_{H'}(u)+\sqrt \epsilon  t,\]
completing the proof.
\hfill$\Box$

\subsection{Monochromatic components}

In \cite{conj-schelp}, Schelp conjectured that if $G$ is a graph on $3n-1$ vertices with minimum degree at least $3|V(G)|/4$, then $G\rightarrow(P_{2n},P_{2n})$ provided $n$ sufficiently large.
Gy\'{a}rf\'{a}s and S\'{a}rk\"{o}zy \cite{gyar-sar-2012} and independently  Benevides, {\L}uczak, Scott, Skokan and White \cite{ben} confirmed this conjecture asymptotically. Recently, Balogh, Kostochka, Lavrov and Liu \cite{balogh-2019} obtain the following result and thus confirm this conjecture thoroughly.
\begin{lemma}[Balogh et al. \cite{balogh-2019}]\label{balogh-conj}
Let $G$ be a graph on $3n -1$ vertices with minimum degree at least $(3|V(G)|-1)/4$. If $n$ is sufficiently large, then $G\rightarrow(P_{2n},P_{2n})$.
\end{lemma}

Erd\H{o}s and Rado remarked that any 2-colored complete graph contains a monochromatic spanning tree, see \cite{gy}. We will apply the following result due to Gy\'{a}rf\'{a}s and S\'{a}rk\"{o}zy \cite{gyar-sar-2012} to get a large monochromatic component for every $2$-coloring of the edges of graph $G$ with large minimum degree.

\begin{lemma}[Gy\'{a}rf\'{a}s and S\'{a}rk\"{o}zy \cite{gyar-sar-2012}]\label{gya-lemma}
For any $2$-coloring of edges of a graph $G$ with minimum degree $\delta(G)\geq\frac{3|V(G)|}{4}$,
there is a monochromatic component of order larger than $\delta(G)$. This estimate is sharp.
\end{lemma}

\section{Proof of Theorem \ref{main theorem}}\label{chap3}

%We focus on the upper bound and assume that $n$ is odd and sufficiently large.
Let $N=(3+\gamma)n$, where $\gamma>0$ is a sufficiently small real number and $n$ is a sufficiently large odd integer.
Consider a 3-edge coloring of $K_N$ on vertex set $V$, and let $G_i$ ($i=1,2,3$) be the graph induced by the edges in the $i$th color. We shall show that either $G_i$ contains a copy of $\mathcal{H}$ for some $i=1,2$, or $G_3$ contains an odd cycle $C_n$.
On the contrary, we suppose that $G_i$ contains no $\mathcal{H}$  for $i=1,2$, and $G_3$ contains no odd cycle $C_n$. We aim to find a contradiction.

We write  $a \ll b$ if $a$ is much smaller than $b$.
Let $\gamma>0$ and $\Delta\geq1$ be given. We apply Lemma \ref{Embedding Lemma} to $G_1$ (or $G_2$)
 with $d=1/4$ and $\Delta$ to get $\epsilon_1$. We set $\eta,\epsilon$,  $\delta$ and $\beta$ such that
\begin{align}\label{eta-ep}
\eta=\frac{\gamma}{15}, \;\;\delta=\min\left\{10^{ -16},\; \frac{\eta^{2}}{100}\right\}
\;\;\epsilon  =\min\left\{\frac{\delta}{1296}, \epsilon_1 \right\}\;\;\text{and}\;\;0<\beta\ll \epsilon.
\end{align}
Moreover, set
\begin{align}\label{constant-1}
t_0= \max \left\{\frac1\epsilon,\;\exp\left(\frac{3(1+\eta)}{\delta^{16}}\right)\right\}.
\end{align}

We apply the regularity lemma (Lemma \ref{regular lemma})
to $G_1$, $G_2$ and $G_3$
with $\epsilon, t_0$ to obtain a $T_0=T_0(\epsilon,t_0)$ such that there exists a partition of the vertex set $V$ into $t+1$ classes $V=V_{0}\cup V_{1}\cup\dots \cup V_t$ satisfying $t_0\le t\le T_0$ and
(1) $\left| {{V_0}} \right| < \epsilon n$, $\left| {{V_1}} \right| = \left| {{V_2}} \right| =  \ldots  = \left| {{V_t}} \right|$;
(2) all but at most $\epsilon {t\choose 2}$ pairs $(V_i,V_j)$, $1 \le i \neq j \le t$, are $\epsilon$-regular for $G_1, G_2$ and $G_3$.
We construct the reduced graph $H$ with vertex set $\{v_1, v_2, \dots, v_t\}$ and the edge set formed by pairs $\{v_i,v_j\}$ for which $(V_i,V_j)$ is $\epsilon$-regular
with respect to $G_1$, $G_2$ and $G_3$.
In which $v_i$ and $v_j$ are non-adjacent in
$H$ if the pairs $(V_i, V_j)$ is not  $\epsilon$-regular for some $G_i$.
Thus we obtain a bijection $f:{v_i} \to {V_i}$ between the vertices of $H$ and the clusters of the partition.

Clearly, $e(H)\geq (1-\epsilon)$$t\choose{2}$ from Lemma \ref{regular lemma}.
Accordingly, we assign color $i$ ($i\in[3]$) to an edge of $H$ if and only if $i$ is the minimum integer for which \[d_{G_i}(V_i,V_j)\geq 1/3.\]
Let $H_i$ ($F_i$) be the spanning subgraph of $H$ ($F$) induced by all edges that have received color $i$.

\medskip
\noindent
{\textbf{Overview of the remaining proof}:} The remaining part of the proof is straightforward but rather rich in technical details, so we shall briefly outline it first.
From the assumption that $G_3$ contains no odd cycle $C_n$, we will show that the reduced subgraph $H_3$ contains no odd cycle of length at least $\frac{{(1 + 0.1\eta )t}}{{3(1 + \eta )}}$. Since $G_i$ contains no $\mathcal{H}$  for $i=1$ or $2$, we shall get that $H_i$  contains no connected matching on more than $(\frac{1}{3}-0.2\eta)t$ vertices  for $i=1$ or $2$. Then, we can easily get that $e(H_3)\geq e(H)/3\geq(\frac{1}{6}-\epsilon/3)t^2$.
Now, since $H_3$ contains no odd cycle of length at least $\frac{{(1 + 0.1\eta )t}}{{3(1 + \eta )}}$, we apply Lemma \ref{luc-3} to graph $H_3$ with $\alpha=\frac{{(1 + 0.1\eta )}}{{3(1 + \eta )}}$ to deduce that $H_3$ contains subgraphs ${H_3'}$ which is bipartite  and $H_3''$ such that the desired properties hold.
Let $A$ and $B$ be the color classes of the bipartition of ${H_3'}$, and let $X=V({H_3''})$.
We will show that $|X|<(1/3-\eta/5)t$ and $\max\{|A|,|B|\}<({\frac{1}{2} - \sqrt \delta/2  } )t.$

Combining with the above desired properties and Lemma \ref{density-degree-upper}, we can  deduce that the reduced graph $H$ contains a subgraph $F$ on at least $(1-\frac{3}{2}\sqrt{\delta})t$ vertices such that $\delta(F)>v(F)  - \sqrt \delta t/2$ and $\deg_{F_3}(u)<(\frac{1}{2} + \sqrt \delta  )t$ for any vertex $u\in V(F)$.
Since $F\subset H$ and $H_i$  contains no connected matching on more than $(\frac{1}{3}-0.2\eta)t$ vertices  for $i=1$ or $2$, we conclude that
$F_i$ contains no connected matching on more than  $(\frac{1}{3} - 3\sqrt \delta  )v({F})$ vertices for $i=1, 2$.
Moreover, by Lemma \ref{su-lemma}, the order of the  largest component of
$F_i ~(i=1, 2)$ must be less than $\frac{1}{3}v(F)$.
By noting $|X|<(1/3-\eta/5)t$
and Fact  \ref{density-degree}, we obtain a
vertex set $B''\subseteq B\cap V(F)$ such that each vertex of $B''$ is adjacent to at least $\max\{v(F)/3, (1-\sqrt{\delta}/4)|B''|\}$.
Now we apply Lemma \ref{gya-lemma} to $H[B'']$ to conclude that $H[B'']$ contains a  monochromatic component in color 1 or 2 of order larger than $\frac{1}{3}v(F)$.
Then $F_i$ (i = 1, 2) contains a monochromatic component in
color $1$ or $2$ of  order at least $\frac{1}{3}v(F)$, which will leads to a contradiction.\hfill$\Box$

\medskip\noindent
{\textbf{Details of the remaining proof}:} We will have the following claims at first.

\begin{claim}\label{long path}
$H_3$ contains no odd cycle of length at least $\frac{{(1 + 0.1\eta )t}}{{3(1 + \eta )}}$.
\end{claim}
\noindent{\bf Proof of Claim \ref{long path}.}
On the contrary, let $C$ be an odd cycle of length $s\ge\frac{{(1 + 0.1\eta )t}}{{3(1 + \eta )}}$. Without loss of generality, assume that its vertex set is $[s]$.
Thus we have that for $i=1,2,\dots,s$,
$(V_i,V_{i+1})$ is $\epsilon$-regular and $d_{G_3}(V_i,V_{i+1})\geq 1/3$, where the indices take modula on $s$.
By Fact \ref{localregular}, for odd $i=1,3,\dots,s-2$,
we can take $V_i'\subseteq V_i$
with $|V_i'|\ge (1-\epsilon)|V_i|$ and $|V_{i+1}'|\ge (1-\epsilon)|V_{i+1}|$ such that $(V'_i,V'_{i+1})$ is $(2\epsilon,1/3-2\epsilon)$-super-regular.
From the property of $\epsilon$-regularity pairs, we can find an odd cycle $u_1u_2\dots u_{s}u_1$ such that $u_i\in V_i'$ for $i=1,2,\dots,s-1$. We apply the techniques used by {\L}uczak \cite{lucazk-1999} and Lemma \ref{long-path-lemma} to show that the $G_3$ contains a $C_n$.

Let $m=(1-\epsilon)^2\frac{N}{t}$, $\epsilon'=2\epsilon$ and $\beta_0=1-6\epsilon.$
Thus we have that  every pair of vertices $v_i\in V_i'$, $v_{i+1}\in V_{i+1}'$ satisfies $\deg_{G_3}(v_i,V_{i+1}')\ge (1/3-2\epsilon)|V_{i+1}'|>\beta_0 m/5$, and $\deg_{G_3}(v_{i+1},V_i')> \beta_0 m/5$.
Therefore, by Lemma \ref{long-path-lemma}, for every $\ell$, $1\le \ell \le m-5\epsilon m/\beta_0$, and for every pair of vertices $u_i\in V_i'$ and $u_{i+1}\in V_{i+1}'$ and odd $i=1,3,\dots,s-2$, $G_3$ contains a path of length $2\ell+1$ connecting $u_i$ and $u_{i+1}$.
Thus, there are odd cycles of all lengths from $s$ to $(s-1)(m-5\epsilon m/\beta_0)$.
Since
\begin{align*}
(s-1)(m-5\epsilon m/\beta_0)&\ge\left(\frac{(1+0.1\eta)t}{3(1+\eta)}-1\right)
\left(1-\frac{5\epsilon} {1-6\epsilon}\right)\left(1-\epsilon\right)^2\frac{N}{t}
\\&\overset{(\ref{constant-1})}{>}\frac{t}{3(1+\eta)}\left(1-10\epsilon \right)\frac{\left(3+\gamma\right)n}{t}\overset{(\ref{eta-ep})}{\ge}\frac{{{\rm{1 + 5}}\eta }}{{{\rm{1 + }}\eta }}\left(1 - 10\epsilon \right)n,
\end{align*}
which is at least $n$ again by noting (\ref{eta-ep}).
So $G_3$ contains an odd cycle $C_n$ as desired.
\hfill$\Box$

\medskip

A connected matching in a graph $G$ is a matching $M$ such that all edges of $M$ are in the same connected component of $G$.
By the assumption of $G_i$ contains no $\mathcal{H}$  for $i=1,2$,
we have the following claim.
\begin{claim}\label{cla-match}
For $i=1,2$, $H_i$ contains no connected matching on more than $(\frac{1}{3}-0.2\eta)t$ vertices.
\end{claim}

\noindent{\bf Proof of Claim \ref{cla-match}}.
On the contrary, without loss of generality, suppose that $H_1$ contains a connected matching  $M$ on at least $(\frac{1}{3}-0.2\eta)t$ vertices that is contained in a tree $T\subset H_1$.
Suppose that the vertex set of $T$ is $\{x_1,\dots,x_l,x_{l+1},\dots, x_{2l}, x_{2l+1},\dots, x_{2l+l'}\}$, and the matching $M$ has edge set
$E_M=\{x_ix_{l+i}: i=1,\dots, l\}$.

We will prove that there exists a copy of $\mathcal{H}$ in $G_1$, contradicting the assumption that $G_1$ contains no $\mathcal{H}$.
Since the proof is similar as in \cite[Theorem 1.3]{G-G-M}, we only give a sketch of the proof as follows.

Firstly,  we shall apply Fact \ref{localregular} to get a  subgraph $G_P$ of $G_1$
with classes $$A_1,\dots,A_{l}, A_{l+1},\dots,A_{2l}, A_{2l+1},\dots, A_{2l+l'}$$ which  corresponds to the vertices
$x_1,\dots,x_l,x_{l+1},\dots, x_{2l}, x_{2l+1},\dots, x_{2l+l'}$, and each of those sets has
 size at least $(1-2\epsilon)N/t$ and  the bipartite graphs induced by $A_i$ and $A_{l+i}$
are $(2\epsilon, 1/3-\epsilon)$-super-regular for $i\in[l]$ and the bipartite graphs induced by all the other pairs are $(2\epsilon, 1/3-\epsilon)$-regular.
It is clear that these dense super-regular pairs covering $(1+o(1))n$ vertices.

Secondly, we partition the vertices of $\mathcal{H}$ and, since $\mathcal{H}=(W,E_\mathcal{H})$
has small bandwidth,
we can apply  Lemma \ref{balanced} to obtain a partition of $W$
which will be composed of clusters $$W_1,\dots,W_l,W_{l+1}\dots,W_{2l}, W_{2l+1},\dots,W_{2l+l'}.$$
Fix $1\leq j\leq 2l+l'$,  we define $U_j$ as the set of vertices of $W_j$ with neighbors in some $W_k$ with $j\neq k$
and $\left\{ {x_jx_k} \right\} \notin M$.
Define the set $U_j'=N_H(U)\cap (W_j\setminus U)$, where $U = \bigcup\nolimits_{I = 1}^{2l + l'} {{U_i}}$.

Then, we  can verify that all of the four conditions of Definition \ref{def} hold, i.e.,

\medskip
(1) $|W_i|\leq |A_i|$ for $i\in[2l+l']$.

\medskip
(2) $xy\in E_\mathcal{H}$ for $x\in W_i, y\in W_j$ implies $x_ix_j\in E_T$ for all distinct $i, j\in [2l+l']$.

\medskip
(3) $|U_i|\leq \epsilon |A_i|$ for $i\in[2l+l']$.

\medskip
(4) $|U'_i|, |U'_j|\leq \epsilon \min\{|A_i|, |A_j|: x_ix_j\in E_{M}\}$.

\medskip
Thus, the partition $\{W_1,\dots, W_{2l+l'}\}$ of $W$
is $(2\epsilon, T , M)$-compatible with $\{A_1,\dots, A_{2l+l'}\}$, which is a partition of $V(G_P)$.

Finally, we can find a copy of $\mathcal{H}$ in $G_P$ by the Embedding Lemma (Lemma \ref{Embedding Lemma}), completing the proof.
\hfill$\Box$

\begin{claim}\label{eh3}
$e(H_3)\geq e(H)/3\geq(\frac{1}{6}-\epsilon/3)t^2$.
\end{claim}

\noindent{\bf Proof of Claim \ref{eh3}.}
On the contrary, suppose that $e(H_3)< e(H)/3$. Thus, without loss of generality, suppose that
\begin{align}\label{eh-3}
e(H_1)>e(H)/3\geq
\frac{1}{3}(1-\epsilon){t\choose{2}}\geq\left(\frac{1}{6}-\epsilon/3\right)t^2.
\end{align}
However, Claim \ref{cla-match} implies that $H_1$ contains no  path with more than $(\frac{1}{3}-0.2\eta)t+1$ vertices, it follows by Lemma \ref{gallai} that
\[
e({H_1}) \leq \frac{(\frac{1}{3}-0.2\eta)t-1}{2}\cdot t< \left(\frac{1}{6} - 0.1\eta \right){t^2}.
\]
This contradicts (\ref{eh-3}) by noting (\ref{eta-ep}).
\hfill$\Box$

\medskip

By Claim \ref{long path}, $H_3$ contains no odd cycle of length at least $\frac{{(1 + 0.1\eta )t}}{{3(1 + \eta )}}$. Now we apply Lemma \ref{luc-3} to graph $H_3$ with $\alpha=\frac{{(1 + 0.1\eta )}}{{3(1 + \eta )}}$ to deduce that $H_3$ contains subgraphs ${H_3'}$ and $H_3''$ such that

\smallskip
($i$) $V({H_3'}) \cup V({H_3''})=V(H_3)$, $V({H_3'})\cap V({H_3''})=\emptyset$ and each of the sets $V({H_3'})$ and $V({H_3''})$  is either empty or contains at least $\alpha \delta t/2$ vertices;

\smallskip

($ii$) ${H_3'}$ is bipartite;

\smallskip

($iii$) ${H_3''}$ contains no more than $\alpha t\left| {V({H_3''})} \right|/2$ edges;

\smallskip

($iv$) all except no more than $\delta {t^2}$ edges of $H_3$ belong to either ${H_3'}$ or ${H_3''}$.

\smallskip

 Let $A$ and $B$ be the sets of the bipartition of ${H_3'}$, and let $X=V({H_3''})$.
 From Claim \ref{eh3} and the property of $H_3$ (see ($iv$)), we have
 \begin{align}\label{eh3-lower}
 e({H_3'})+e({H_3''})&\geq e(H_3)-\delta t^2\geq\left(\frac{1}{6}-\epsilon/3 \right)t^2-\delta t^2 \overset{(\ref{eta-ep})}{>}\left(\frac{1-7\delta}{6}\right)t^2.
 \end{align}
%where the last inequality holds as $\epsilon<\frac{\delta}{2}$ from (\ref{eta-ep}).

\begin{claim}\label{claim123}
$|X|<(1/3-\eta/5)t$.
\end{claim}
\noindent{\bf Proof of Claim \ref{claim123}.}  Let us put $|X|=\lambda t$ and $|A\cup B|=(1-\lambda)t$. By Lemma \ref{luc-3} ($iii$),
\begin{align}\label{eh3-upper}
e({H_3'})+e({H_3''}) &{\le}  e({H_3'})+ \frac{{\alpha t\left| {V({H_3''})} \right|}}{2} \le\frac{|V(H_3')|^2}{4}+\frac{{\alpha t\left| {V({H_3''})} \right|}}{2}\nonumber \\&\le\frac{{{{\left( {1 - \lambda } \right)}^2}}}{4}{t^2}+\frac{{\lambda \left( {1 + 0.1\eta } \right)}}{{6(1 + \eta )}}{t^2},
\end{align}
which together with (\ref{eh3-lower}) yield that
\[\frac{{{{\left( {1 - \lambda } \right)}^2}}}{4}+\frac{{\lambda \left( {1-0.8\eta } \right)}}{6} >\frac{{{{\left( {1 - \lambda } \right)}^2}}}{4}+\frac{{\lambda \left( {1 + 0.1\eta } \right)}}{{6(1 + \eta )}} \geq\frac{1}{6}-\frac{7\delta}{6},
\]
from which we obtain that $3\lambda^2-(4+1.6\eta)\lambda+(1+14\delta)>0.$
Since $\lambda\leq 1$,  it follows that
\begin{align*}
\lambda<\frac{1}{6}\left(4+1.6\eta-\sqrt{(4+1.6\eta)^2-12(1+14\delta)}\right)
<\frac{1}{3}\left(2+0.8\eta-\sqrt{1+3\eta}\right)<\frac13-\frac\eta5
\end{align*}
provided $\eta$ is sufficiently small and by noting (\ref{eta-ep}).
\hfill$\Box$

\begin{claim}\label{claim12}
$\max\{|A|,|B|\}<({\frac{1}{2} - \sqrt \delta/2  } )t.$
\end{claim}

\noindent{\bf Proof of Claim \ref{claim12}.}
On the contrary, suppose that $|A|\geq ( {\frac{1}{2} - \sqrt \delta /2 } )t$ without loss of generality.
Let $H[A]$ be the subgraph of $H$ induced by $A$.
Note that $H[A]$ only contains edges in color 1 or color 2 as $H_3'$ is bipartite.
Note also that $e(H[A])\geq{|A|\choose{2}}-\epsilon{t\choose{2}}.$
By Fact \ref{density-degree}, $H[A]$ contains a subgraph $H[A']$  such that
$|A'|\ge({\frac{1}{2} - \sqrt \delta/2  } )t-\sqrt{\epsilon}t>( \frac{1}{2}-\sqrt{\delta} )t$ and
\begin{align*}
\delta(H[A']) &\ge \left| {A'} \right| - 2\sqrt \epsilon  t > \left| {A'} \right| - 6\sqrt \epsilon  | A'|\overset{(\ref{eta-ep})}{>}(1 - \sqrt \delta)| {A'}|,
\end{align*}
the second inequality due to $t<\frac{|A'|}{(1/2-\sqrt \delta)}$ and $\delta$ is sufficiently small.
Now we apply Lemma \ref{balogh-conj} with graph $H[A']$ and $n=\frac{(1/2-\sqrt{\delta})t+1}{3}$ to conclude that $H[A']$ contains a monochromatic path
${P_{\frac{{(1- 2\sqrt \delta )t }}{3}}}$ in color 1 or color 2.
As $\frac{(1 - 2\sqrt \delta )t}{3}>(\frac{1}{3}-0.2\eta)t+1$ from (\ref{eta-ep}), we have that either $G_1$ or $G_2$ must contain a copy of $\mathcal{H}$ by Claim \ref{cla-match}, a contradiction.
\hfill$\Box$

\medskip

For graphs $P$ and $Q$, we use $P\cup Q$ to denote the graph defined on $V(P)\cup V(Q)$ whose edge set is $E(P)\cup E(Q)$.
From Claims \ref{claim123}-\ref{claim12}  noting that $V({H_3'}) \cup V({H_3''})=V(H_3)=V(H)$ and $V({H_3'})\cap V({H_3''})=\emptyset$, we have that for any vertex $u\in V(H)$,
\begin{align}\label{degree-upper}
\deg_{H_3'\cup H_3''}(u)< \left( {\frac{1}{2} - \frac{\sqrt \delta}2} \right)t.
\end{align}

For an edge-colored graph  $F$, we use $F_i$ to denote the subgraph induced by edges in color $i$ in $F$.

\begin{claim}\label{cla6}
$H$ has a subgraph $F$ with at least $(1-\frac{3}{2}\sqrt{\delta})t$ vertices such that $\delta(F)>v(F)  - \sqrt \delta t/2$ and $\deg_{F_3}(u)<(\frac{1}{2} + \sqrt \delta  )t$ for any vertex $u\in V(F)$.
\end{claim}
\noindent{\bf Proof of Claim \ref{cla6}.}
We know that all but at most $\delta t^2$
edges of $H_3$ are contained in $H_3' \cup H_3''$.
By Lemma \ref{density-degree-upper}, there exists an induced subgraph $H_3^0$ of $H_3$ such that  $v(H_3^0)\geq v(H)-\sqrt{2\delta}t$
and
\begin{align}\label{q-degree}
\deg_{H_3^0}(u)<\deg_{{H_3'}\cup{H_3''}}(u)+\sqrt{2\delta}t\overset{(\ref{degree-upper})}{<}\left( {\frac{1}{2} - \sqrt \delta/2} \right)t+\sqrt{2\delta}t<\left(\frac{1}{2}+\sqrt{\delta}\right)t
\end{align}
for any vertex $u\in V(H_3^0)$.

Note that $e(H[V(H_3^0)])\geq {v(H_3^0)\choose{2}}-\epsilon{t\choose{2}}$,
it follows from Fact \ref{density-degree} that $H[V(H_3^0)]$ contains a subgraph $F$ such that
\begin{align}\label{F-degree}
v(F)\geq v(H_3^0)-\sqrt{\epsilon}t\geq v(H)-\sqrt{2\delta}t-\sqrt{\epsilon}t\overset{(\ref{eta-ep})}{>}\left(1-\frac{3}{2}\sqrt{\delta}\right)t
\end{align}
and
\begin{align}\label{F-degree-min}
\delta(F)\geq v(H_3^0)-2\sqrt{\epsilon}t>v(F)-\sqrt \delta t/2.
\end{align}
Since $F_3\subset H_3^0$, by noting (\ref{q-degree}), we obtain that for any $u\in V(F)$,
\begin{align}\label{f3-upper-degree}
\deg_{F_3}(u)\leq \deg_{H_3^0}(u)<\left(\frac{1}{2}+\sqrt{ \delta}\right)t,
\end{align}
completing the proof.
\hfill $\Box$

\medskip
According to  (\ref{F-degree-min}) and (\ref{f3-upper-degree}),
we obtain that for any vertex $u\in V(F)$,
\begin{align}\label{degree-F12}
\deg_{F_1\cup F_2}(u)&={\deg_{{F_1}}}(u) + {\deg_{{F_2}}}(u) \geq \delta(F)-d_{F_3}(u)\nonumber\geq v(F)  - \sqrt \delta t/2-\left(\frac{1}{2}+\sqrt{ \delta}\right)t\nonumber \\&\overset{(\ref{F-degree})}{>} \left( {1 - \frac{{\frac{1}{2}+\frac{3}{2}\sqrt \delta }}{{1 - \frac{3}{2}\sqrt \delta  }}} \right)v(F)> \left( \frac{1}{2} -3\sqrt {\delta} \right) v(F).
\end{align}

 Recall that a connected matching in a graph $G$ is a matching $M$ such that all edges of $M$ are in the same connected component of $G$.
\begin{claim}\label{cla7}
For $i=1,2$, $F_i$ contains no connected matching on more than  $(\frac{1}{3} - 3\sqrt \delta  )v({F})$ vertices.
\end{claim}
\noindent{\bf Proof of Claim \ref{cla7}.}
Note that since $F\subset H$, by Claim  \ref{cla-match}, $F_i$ contains no connected
matching on  more than $(\frac{1}{3}-0.2\eta)t$ vertices for $i=1,2$.
Since
\[
\left({\frac{{\rm{1}}}{{\rm{3}}}{\rm{ - 3}}\sqrt \delta  } \right)v(F) \overset{(\ref{F-degree})}{>} \left( {\frac{{\rm{1}}}{{\rm{3}}}{\rm{ - 3}}\sqrt \delta  } \right)\left( {1 - \frac{{3\sqrt \delta  }}{2}} \right)t > \left( {\frac{{\rm{1}}}{{\rm{3}}} - \frac{{7\sqrt \delta  }}{2}} \right)t \overset{(\ref{eta-ep})}{> } \left(\frac{{\rm{1}}}{{\rm{3}}} - 0.2\eta \right)t.
\]
 Thus  $F_i$  ($i=1,2$) contains no connected matching on more than $(\frac{1}{3} - 3\sqrt \delta  )v({F})$ vertices.
\hfill $\Box$

 \begin{claim}\label{lag-cp-F12}
For $i=1,2$, the largest component of $F_i$ has order less than $\frac{1}{3}v(F)$.
\end{claim}
\noindent{\bf Proof of Claim \ref{lag-cp-F12}.}
Let $R_1,\dots, R_r$ be the components of $F_1$ and let $B_1,\dots,B_b$ be the components of $F_2$. Without loss of generality, suppose that $|V(R_i)|\geq |V(R_{i+1})|$ for all $1\leq i\leq r-1$, and $|V(B_j)\geq |V(B_{j+1})|$ for all
$1\leq j\leq b-1$.
Let $r'$, $0\leq r' \leq r $, be the maximum integer such that $|V(R_{r'})|\geq (\frac{1}{3} - 3\sqrt \delta)v(F)$.
Similarly, let $b'$, $0\leq b' \leq b$, be the maximum integer such that
$|V(R_{b'})|\geq(\frac{1}{3} - 3\sqrt \delta)v(F)$. We aim to show that $r'=b'=0$.

By Claim \ref{cla7}, $F_1$ and hence each $R_i$ contains no connected matching on more than $(\frac{1}{3}-3\sqrt \delta  )v({F})$ vertices.
Thus Lemma \ref{su-lemma} implies that each $R_i$ has a partition
$S_{R_i}\cup Q_{R_i}\cup I_{R_i}$ satisfying
\begin{align}\label{su-ineq-1}
\left|{{Q_{{R_i}}}} \right| + 2\left| {{S_{{R_i}}}} \right| = \min \left\{ {\left| {V({R_i})} \right|,\left( {\frac{1}{3} - 3\sqrt \delta  } \right)v(F) - 1} \right\}.
\end{align}

\begin{prop}\label{cla-empty}
For $i>r'$,  $I_{R_i}=\emptyset$ and $S_{R_i}=\emptyset$; similarly, for $j>b'$,  $I_{B_j}=\emptyset$ and $S_{B_j}=\emptyset$.
\end{prop}
{\bf Proof.} For  $i>r'$,  $\left| {V({R_i})} \right| \le ( {\frac{1}{3} - 3\sqrt \delta  } )v(F)-1$, implying that $\left| {{Q_{{R_i}}}} \right| + 2\left| {{S_{{R_i}}}} \right| = \left| {V({R_i})} \right|$ and $I_{R_i}=\emptyset$ by  (\ref{su-ineq-1}).
Thus $\left| {{Q_{{R_i}}}} \right| + \left| {{S_{{R_i}}}} \right| = \left| {V({R_i})} \right|$,
and $S_{R_i}=\emptyset$ follows. The second assertion is similar.
\hfill$\Box$

\medskip
According to Lemma \ref{su-lemma} and (\ref{su-ineq-1}), we obtain that for $1\le i\le r$,
\begin{align}\label{whit-req1}
|S_{R_i}|<\frac{1}{2}\left(\frac{1}{3}- 3\sqrt \delta \right)v(F),
\;\;|Q_{R_i}|+2|S_{R_i}|<\left(\frac{1}{3}- 3\sqrt \delta  \right)v(F),
\end{align}
and
\begin{align}\label{whit-req2}
\nonumber	|I_{R_i}|=|V(R_i)\setminus (Q_{R_i}\cup S_{R_i})|&=|V(R_i)|+|S_{R_i}|-(|Q_{R_i}|+2|S_{R_i}|)
\\&>|V(R_i)|+|S_{R_i}|-\left(\frac{1}{3}- 3\sqrt \delta  \right)v(F).
\end{align}
Note that for any $B_j$,  similar sets $S_{B_j}$, $I_{B_j}$ and $Q_{B_j}$ can be defined, with analogues of the above bounds. In particular,
$| {{I_{B_j}}}|\ge { | {V({B_j})} |+{| {{S_{{B_j}}}} |  -  (\frac{1}{3} - 3\sqrt \delta  )v(F)} }$.

By Lemma \ref{su-lemma} ($iv$), for $1\leq i\leq r$, any vertex $u\in I_{R_i}$ satisfies that $\deg_{F_1}(u)\le\frac{1}{2}(\frac{1}{3} - 3\sqrt \delta  )v(F)$.
Similarly, for $1\leq j\leq b$, any vertex  $u\in I_{B_j}$ satisfies that $\deg_{F_2}(u)\le\frac{1}{2}(\frac{1}{3} - 3\sqrt \delta  )v(F)$.

Note that for $1\leq i\leq r$, each vertex in $Q_{R_i}$ has at most one neighbor in $I_{R_i}$ in $F_1$ by Lemma \ref{su-lemma} ($iii$).
Hence, any vertex $u\in Q_{R_i}$ satisfies that $\deg_{F_1}(u)\le|(Q_{R_i}\cup S_{R_i})\setminus\{u\}|+1$, which is less than $(\frac{1}{3} - 3\sqrt \delta)v(F)$ by noting (\ref{whit-req1}).
Similarly, for  $1\leq j\leq b$, any vertex $u\in Q_{B_j}$ satisfies that $\deg_{F_2}(u)<(\frac{1}{3} - 3\sqrt \delta  )v(F)$.

\begin{prop}\label{white-claim2}
 $I_{R_i}\cap(I_{B_j}\cup Q_{B_j})= \emptyset $.
\end{prop}
{\bf Proof.} Indeed, if there is a vertex $u\in I_{R_i}\cap(I_{B_j}\cup Q_{B_j})$, then from the above observation,
$$\deg_{F_1\cup F_2}(u)=\deg_{F_1}(u)+\deg_{F_2}(u)<\frac{1}{2}\left(\frac{1}{3}-3\sqrt \delta\right)v(F)+\left(\frac{1}{3}-3\sqrt \delta\right)v(F),$$
which is less than $(\frac{1}{2}-4\sqrt\delta)v(F)$, contradicting (\ref{degree-F12}).
\hfill$\Box$

\medskip

Therefore, for $1\leq i\leq r$, the set $I_{R_i}\subseteq \bigcup\nolimits_{1 \le j \le b} {{S_{{B_j}}}}$
due to Proposition \ref {white-claim2}.
Since $S_{B_j}=\emptyset$ for $j>b'$ by Proposition \ref{cla-empty}, we have that $I_{R_i}\subseteq \bigcup\nolimits_{1 \le j \le b'} {{S_{{B_j}}}}$.
Similarly, for $1\leq j\leq b$, $I_{B_j}\subseteq\bigcup\nolimits_{1 \le i \le r'} {{S_{{R_i}}}} $. Consequently, again by Proposition \ref{cla-empty}, $I_{R_i}=\emptyset$ for $i>r'$ and $I_{B_j}=\emptyset$ for  $j>b'$, we obtain that
\begin{align}\label{B_r_R}
\sum\limits_{i= 1}^{r'} {\left| {{I_{{R_i}}}} \right|}=\sum\limits_{i = 1}^{r} {\left| {{I_{{R_i}}}} \right|}  \le \sum\limits_{j = 1}^{b'} {\left| {{S_{{B_j}}}} \right|},\;\;\text{and}\;\;\sum\limits_{j= 1}^{b'} {\left| {{I_{{B_j}}}} \right|}  \le \sum\limits_{i = 1}^{r'} {\left| {{S_{{R_i}}}} \right|},
\end{align}
which together with  (\ref{whit-req2}) yield that
\begin{align}\label{Bj-Ri}
\nonumber0 &\ge \sum\limits_{j = 1}^{b'} {\left| {{I_{B_j}}} \right|}  - \sum\limits_{i = 1}^{r'} {\left| {{S_{{R_i}}}} \right|}		\nonumber\\& \ge\sum\limits_{j = 1}^{b'} {\left(\left| {V({B_j})} \right|+ {\left| {{S_{{B_j}}}} \right|   -  \left(\frac{1}{3} - 3\sqrt \delta  \right)t} \right)}  - \sum\limits_{i = 1}^{r'} {\left| {{S_{{R_i}}}} \right|}		\nonumber\\& \overset{(\ref{B_r_R})}{\ge} \sum\limits_{j = 1}^{b'} {\left( {\left| {V({B_j})} \right| -  \left(\frac{1}{3} - 3\sqrt \delta  \right)t } \right)}  + \sum\limits_{i = 1}^{r'} {\left| {{I_{{R_i}}}} \right|}  - \sum\limits_{i = 1}^{r'} {\left| {{S_{{R_i}}}} \right|}
	\nonumber\\&\overset{(\ref{whit-req2})}{>} \sum\limits_{j = 1}^{b'} {\left( {\left| {V({B_j})} \right| -  {\left(\frac{1}{3} - 3\sqrt \delta  \right)v(F)  } } \right)}  + \sum\limits_{i = 1}^{r'} {\left( {\left| {V({R_i})} \right| -  {\left(\frac{1}{3} - 3\sqrt \delta  \right)v(F) } } \right)},
\end{align}
which implies that $r'=b'=0$ since $\left| {V({R_{i}})} \right| \ge (\frac{1}{3} - 3\sqrt \delta  )v(F)$ and $\left| {V({B_{j}})} \right| \ge (\frac{1}{3} - 3\sqrt \delta  )v(F)$ for $1\leq i\leq r'$ and $1\leq j\leq b'$.
This completes Claim \ref{lag-cp-F12}.
\hfill$\Box$

\medskip

By noting Claim \ref{claim123} and $A\cup B\cup X=V(H_3)=V(H)$ and $F\subseteq H$, we obtain that
\begin{align*}
\left|{\left( A \cup B \right) \cap V({F})} \right|&\ge |V({F})|-|X|> v({F}) - \left(\frac{1}{3} - \frac{1}{5}\eta \right)t\overset{(\ref{F-degree})}{>} v(F) - \frac{{\frac{1}{3} - \frac{1}{5}\eta }}{{1 - \frac{3}{2}\sqrt \delta  }}v(F)
\\&\overset{(\ref{eta-ep})}{>}v(F) - \frac{1}{3}(1 - \sqrt \delta/2  )v(F)= \left( {\frac{2}{3} + \frac{1}{6}\sqrt \delta } \right)v(F).
\end{align*}
Thus, one of $|A \cap V({F})|$ and $|B \cap V(F)|$ must be at least
$({\frac{1}{3} + \frac{1}{12}\sqrt \delta  } )v(F).$
Without loss of generality,  we suppose that ${B'} \subseteq {B} \cap V({F})$
satisfies
\begin{align}\label{v3equ}
\left|B' \right| = \left( {\frac{1}{3} + \frac{1}{12}\sqrt \delta  } \right)v(F).
\end{align}
Note that $e(H[B'])\geq {|B'|\choose{2}}-\epsilon{t\choose{2}}.$
Therefore, by deleting at most $\sqrt{{\epsilon}}t$ vertices from $B'$ we obtain a vertex set $B''\subseteq B'$ such that each vertex of $B''$ is adjacent to at least $|B'|-2\sqrt{\epsilon}t$ in $H[B'']$ from Fact \ref{density-degree}.
We obtain that
\begin{align}\label{last-ineq1}
\delta (H[B'']) \ge \left| {B'} \right| - 2\sqrt{\epsilon}t\overset{(\ref{v3equ}),  (\ref{F-degree})}{>} \left( {\frac{1}{3} + \frac{1}{12}\sqrt \delta  } \right)v(F) - \frac{{2\sqrt {\epsilon } }}{{1 - \frac{3}{2}\sqrt \delta  }}v(F)
\overset{(\ref{eta-ep})}{\geq} \frac{1}{3} v(F).
\end{align}
By noting (\ref{v3equ}),
\begin{equation*}
\delta (H[B'']) >\frac{1}{3}v(F)= \frac{{\left| {B'} \right|}}{{1 + \sqrt \delta/4 }}
\ge\frac{{\left|{B''} \right|}}{{1 + \sqrt \delta/4 }}> \left( {1 - \frac{ \sqrt{\delta}  }{4}} \right)\left| {B''} \right|.
\end{equation*}

Note that all edges of $H[B'']$ are colored only with color $1$ or $2$ since $B''\subseteq B$ and $B$ is one of parts of the bipartite graph $H_3'$.
Now we apply Lemma \ref{gya-lemma} to $H[B'']$ to conclude that $H[B'']$ contains a  monochromatic component in color 1 or 2 of order larger than $\delta (H[B''])$, which is at least $\frac{1}{3}v(F)$ according to (\ref{last-ineq1}).
Since  $H[B''] \subseteq {F}$, it follows that $F$ contains a monochromatic component in color 1 or 2 of order at least $\frac{1}{3}v(F)$.
This contradicts Claim \ref{lag-cp-F12}.

In conclusion, the proof of Theorem \ref{main theorem} is complete.
\hfill $\Box$

\end{spacing}


\begin{thebibliography}{99}

\bibitem{a-b-s-2013}
P. Allen, G. Brightwell and  J. Skokan,   Ramsey-goodness and otherwise, \emph{Combinatorica} 33 (2013),  125--160.

\bibitem{afz}
N. Alon, J. Fox and Y. Zhao,  Efficient arithmetic regularity and removal lemmas for induced bipartite patterns, Discrete Anal. 2019, Paper No. 3, 14 pp.


\bibitem{balogh-2019}
J. Balogh,  A. Kostochka,  M. Lavrov and  X. Liu,   Monochromatic paths and cycles in $2$-edge-colored graphs with large minimum degree, arXiv:1906.02854.

\bibitem{ben}
 F. S. Benevides, T. Luczak, A. Scott, J. Skokan and M. White, Monochromatic cycles in 2-coloured graphs, \emph{Combin. Probab. Comput.}  21  (2012), 57--87.

\bibitem{Benevides-Skokan-2009}
 F. S. Benevides and J. Skokan, The 3-colored Ramsey number of even cycles, \emph{J. Combin. Theory Ser. B} 99 (2009), 690--708.

\bibitem{Bielak}
H. Bielak, Multicolor Ramsey numbers for some paths and cycles,  \emph{Discuss. Math. Graph Theory} 29 (2009), 209--218.

\bibitem{b-e}
J. A. Bondy and P. Erd\H{o}s, Ramsey numbers for cycles in graphs, {\em J. Combin. Theory Ser. B} 14 (1973), 46--54.

\bibitem{J-B}
J. B\"{o}ttcher, {Embedding large graphs--The Bollob\'{a}s-Koml\'{o}s conjecture and beyond}, Ph.D. thesis, Technischen Universit\"{a}t M\"{u}nchen, 2009.

\bibitem{J-P}
J. B\"{o}ttcher, P. Heinig and A. Taraz, Embedding into bipartite graphs, {\em SIAM J. Discrete Math.} 24 (2010), 1215--1233.


\bibitem{julia-klaas}
J. B\"{o}ttcher, K. P. Pruessmann, A. Taraz and A. W\"{u}rfl, Bandwidth, expansion, treewidth, separators and universality for bounded-degree graphs, \emph{European J. Combin.} 31 (2010), 1217--1227.

\bibitem{cly}
X. Chen, Q. Lin and C. You, Ramsey numbers of large books, {\em J. Graph Theory}, to appear. DOI: 10.1002/jgt.22815.


\bibitem{conlon}
D. Conlon, The Ramsey number of books, \emph{Adv. Combin.} 3 (2019), 12pp.

\bibitem{cf}
D. Conlon and J. Fox,  Bounds for graph regularity and removal lemmas, {\em Geom. Funct. Anal.} 22 (2012), 1191--1256.

\bibitem{cfw}
D. Conlon, J. Fox and Y. Wigderson, Ramsey number of books and quasirandomness, {\em Combinatorica}, to appear.



%\bibitem{Davies-2017}
%E. Davies, M. Jenssen, and B. Roberts, Multicolour Ramsey numbers of paths and even cycles, \emph{European J. Combin.}, 63 (2017), 124--133.


%\bibitem{DzNS}
%T. Dzido, A. Nowik and P. Szuca, New lower bound for multicolor Ramsey numbers for even cycles, \emph{Electron. J. Combin.}, 12 (2005), \#N13.

\bibitem{DzFi2}
T. Dzido and R. Fidytek, On some three color Ramsey numbers for paths and cycles, \emph{Discrete Math.} 309 (2009), 4955--4958.




\bibitem{efrs}
P. Erd\H{o}s, R. J. Faudree, C. C. Rousseau and R. H. Schelp, Generalized Ramsey theory for multiple colors, \emph{J. Combin. Theory Ser. B} 20 (1976), 250--264.


\bibitem{Erd-Gallai-1959}
P. Erd\H{o}s and T. Gallai,  On maximal paths and circuits of graphs, \emph{Acta Math. Hungar.} 10 (1959), 337--356.

\bibitem{FLPS}
R. J. Faudree, S. L. Lawrence, T. D. Parsons and R. H. Schelp, Path-Cycle Ramsey numbers, \emph{Discrete Math.} 10 (1974), 269--277.


\bibitem{f-s}
R. J. Faudree and R. H. Schelp, All Ramsey numbers for cycles in graphs, \emph{Discrete Math.} 8 (1974), 313--329.


\bibitem{Faudree-Schelp-1975}
R. J. Faudree and R. H. Schelp,  Path Ramsey numbers in multicolorings, \emph{J. Combin. Theory Ser. B} 19 (1975), 150--160.

\bibitem{f-l-2007}
A. Figaj and T. {\L}uczak, The Ramsey number for a triple of long even cycles, \emph{J. Combin. Theory Ser. B} 97 (2007), 584--596.


\bibitem{Figaj-luczak-2018}
A. Figaj and T. {\L}uczak, The Ramsey numbers for a triple of long cycles, \emph{Combinatorica} 38 (2018), 827--845.

\bibitem{Ferguson-1}
D. G. Ferguson, The Ramsey number of mixed-parity cycles I,  arXiv:1508.07154.

\bibitem{Ferguson-2}
D. G. Ferguson, The Ramsey number of mixed-parity cycles II,  arXiv:1508.07171.

\bibitem{Ferguson-3}
D. G. Ferguson, The Ramsey number of mixed-parity cycles III, arXiv:1508.07176.


\bibitem{flz}
J. Fox, L. M. Lov\'{a}sz and Y. Zhao, On regularity lemmas and their algorithmic applications, {\em Combin. Probab. Comput.} 26 (2017), 481--505.

\bibitem{gerence-gyarfas-1967}
L. Gerencs\'{e}r and A. Gyarf\'{a}s, On Ramsey-type problems, \emph{Ann. Univ. Sci. Budapest E\"{o}tv\"{o}s Sect. Math.} 10 (1967), 167--170.



\bibitem{gy}
A. Gy\'{a}rf\'{a}s, Large monochromatic components in edge colorings of graphs: A survey. In Ramsey Theory: Yesterday, Today and Tomorrow (A. Soifer, ed.), Birkh\"{a}user, pp. 77--96, 2010.


\bibitem{gyarfas-szemerdi-2007}
A. Gy\'{a}rf\'{a}s, M. Ruszink\'{o}, N. S\'{a}rk\"{o}zy and E. Szemer\'{e}di,  Three-color Ramsey numbers for Paths, \emph{Combinatorica} 27 (2007), 35--69.



\bibitem{gyar-sar-2012}
A. Gy\'{a}rf\'{a}s and G. N. S\'{a}rk\"{o}zy, Star versus two stripes Ramsey numbers and a conjecture of Schelp, \emph{Combin. Probab. Comput.} 21 (2012),  179--186.


\bibitem{j-s}
M. Jenssen and J. Skokan, Exact Ramsey numbers of odd cycles via nonlinear optimisation, \emph{Adv. Math.} 376 (2021), 107444, 46 pp.


\bibitem{ks}
J. Koml\'{o}s and M. Simonovits, ``Szemer\'{e}di's regularity lemma and its applications in graph theory,'' Combinatorics, Paul Erd\H{o}s is eighty, Vol. 2 (Keszthely, 1993), 295--352, Bolyai Soc. Math. Stud., 2, J\'{a}nos Bolyai Math. Soc., Budapest, 1996.


\bibitem{kss-2005}
Y. Kohayakawa, M. Simonovits and J. Skokan, The 3-colored Ramsey number of odd cycles, \emph{Electron. Notes Discrete Math.} 19 (2005), 397--402.

\bibitem{kn-su}
C. Knierim and P. Su, Improved bounds on the multicolor Ramsey numbers of paths and even cycles, \emph{Electron. J. Combin.} 26 (2019), \#P1.26.

\bibitem{lp}
Q. Lin and X. Peng, Large book--cycle Ramsey numbers, {\em SIAM J. Discrete Math.} 35 (2021), 532--545.

\bibitem{lucazk-1999}
T. {\L}uczak, $R(C_n, C_n, C_n) \leq (4 + o(1))n$, \emph{J. Combin. Theory Ser. B} 75 (1999), 174--187.


\bibitem{lucss}
T. {\L}uczak, M. Simonovits and J. Skokan, On the multi-colored Ramsey numbers of cycles, \emph{J. Graph Theory} 69 (2012), 169--175.






%\bibitem{li}
%Y. Li, The multicolor Ramsey number of an odd cycle, \emph{J. Graph Theory}, 62 (2009), 324--328.

\bibitem{G-G-M}
 G. Mota, G. N. S\'{a}rk\"{o}zy, M. Schacht and  A. Taraz, {Ramsey number for bipartite graphs with small bandwidth}, {\em European J. Combin.} 48 (2015), 165--176.

\bibitem{nr09}
V. Nikiforov and C. C. Rousseau, Ramsey goodness and beyond, {\em Combinatorica}  29 (2009), 227--262.


\bibitem{Omidi-Raeisi-2011}
G. R. Omidi and G. Raeisi, On multicolor Ramsey number of paths versus cycles, \emph{Electron. J. Combin.}  18 (2011), \#P24.

\bibitem{rs}
V. R\"{o}dl and M. Schacht,  Regularity lemmas for graphs. Fete of combinatorics and computer science, 287--325, Bolyai Soc. Math. Stud., 20, J\'{a}nos Bolyai Math. Soc., Budapest, 2010.

\bibitem{ros}
V. Rosta, On a Ramsey-type problem of J. A. Bondy and P. Erd\H{o}s. I, II, {\em J. Combin. Theory Ser. B} 15 (1973), 105--120.

\bibitem{sark-2016}
G. N. S\'{a}rkozy,  On the multi-colored Ramsey numbers of paths and even cycles, \emph{Electron. J. Combin.} 23 (2016), \#P3.

\bibitem{conj-schelp}
R. H. Schelp, Some Ramsey-Tur\'{a}n type problems and related questions, \emph{Discrete Math.} 312  (2012), 2158--2161.

\bibitem{shao-xu}
Z. Shao, X. Xu, X. Shi and L. Pan, Some three-color Ramsey numbers, $R(P_4, P_5, C_k)$ and $R(P_4, P_6, C_k)$, \emph{European J. Combin.} 30 (2009), 396--403.


\bibitem{regular-lemma}
E. Szemer\'{e}di, Regular partitions of graphs, Probl\`{e}mes combinatoires et
th\'{e}orie des graphes (Colloq. Internat. CNRS, Univ.  Orsay, 1976), 399--401, Colloq. Internat. CNRS, 260, CNRS, Paris, 1978.

\bibitem{sz15}
E. Szemer\'{e}di, Arithmetic progressions, different regularity lemmas and removal lemmas, {\em Commun. Math. Stat.} 3 (2015), 315--328.








\end{thebibliography}
\end{document}